\documentclass[14pt]{article}
\usepackage{mathrsfs}
\usepackage{amsthm}
\usepackage{amssymb}
\usepackage{amsmath}
\usepackage{graphicx}
\usepackage{color}
\usepackage{amsfonts}
\usepackage{float}
\usepackage{cite}
\usepackage[text={140mm,210mm},left=45mm,vmarginratio=1:1]{geometry}
\newtheorem{theorem}{Theorem}[section]

\newtheorem{lemma}[theorem]{Lemma}

\numberwithin{equation}{section}
\normalsize

\begin{document}
\title{\textbf{Phase transition for the SIR model with random transition rates on complete graphs}}

\author{Xiaofeng Xue \thanks{\textbf{E-mail}: xfxue@bjtu.edu.cn \textbf{Address}: School of Science, Beijing Jiaotong University, Beijing 100044, China.}\\ Beijing Jiaotong University}

\date{}
\maketitle

\noindent {\bf Abstract:} In this paper we are concerned with the Susceptible-Infective-Removed model with random transition rates on complete graphs $C_n$ with $n$ vertices. We assign i. i. d. copies of a positive random variable $\xi$ on each vertex as the recovery rates and i. i. d copies of a positive random variable $\rho$ on each edge as the edge infection weights. We assume that a susceptible vertex is infected by an infective one at rate proportional to the edge weight on the edge connecting these two vertices while an infective vertex becomes removed with rate equals the recovery rate on it, then we show that the model performs the following phase transition when at $t=0$ one vertex is infective and others are susceptible. When $\lambda<\lambda_c$, the proportion of vertices which have ever been infective converges to $0$ weakly as $n\rightarrow+\infty$ while when $\lambda>\lambda_c$, there exist $c(\lambda)>0$ and $b(\lambda)>0$ such that for each $n\geq 1$ with probability at least $b(\lambda)$ the proportion of vertices which have ever been infective is at least $c(\lambda)$. Furthermore, we prove that $\lambda_c$ is the inverse of the production of the mean of $\rho$ and the mean of the inverse of $\xi$.

\quad

\noindent {\bf Keywords:} SIR, complete graph, phase transition, random rate.

\section{Introduction}\label{section one}
In this paper we are concerned with the SIR (Susceptible-Infective-Removed) model with random transition rates on complete graphs. First we introduce some notations and definitions to define our model. For each integer $n\geq 1$, we denote by $C_n$ the complete graph with $n$ vertices. We use $0,1,\ldots,n-1$ to denote the $n$ vertices on $C_n$. For $0\leq i\neq j\leq n-1$, we use $e(i,j)$ to denote the edge connecting $i$ and $j$. Note that $e(i,j)=e(j,i)$ for any $i\neq j$. With our notations, $C_n$ is a subgraph of $C_m$ for any integers $n\leq m$. Let $\xi$ be a random variable such that $P(\xi\geq 1)=1$, then we assume that $\{\xi(j)\}_{j=0}^{+\infty}$ are i. i. d. copies of $\xi$. Let $\rho$ be a random variable such that $P(0\leq \rho\leq 1)=1$ and $P(\rho>\epsilon_0)>0$ for some $\epsilon_0>0$, then we assume that $\{\rho(i,j)\}_{0\leq i<j<+\infty}$ are i. i. d copies of $\rho$. We assume that $\{\xi(j)\}_{j=0}^{+\infty}$ and $\{\rho(i,j)\}_{0\leq i<j<+\infty}$ are independent. For $i>j$, we define $\rho(i,j)=\rho(j,i)$. For our model, $\xi(j)$ is the random recovery rate of the vertex $j$ while $\rho(i,j)$ is the random infection weight on the edge $e(i,j)$, for which we assume that $\rho(i,j)=\rho(j,i)$.

The SIR model on $C_n$ is a Markov process with state space $X=\{0,1,-1\}^{C_n}$. That is to say, there is a spin taking a value from $\{0,1,-1\}$ on each vertex of $C_n$. For any $t\geq 0$ and integer $0\leq i\leq n-1$, we denote by $\eta_t$ the state of the process at moment $t$ and $\eta_t(i)$ the value of the spin on $i$ at moment $t$. For any $t\geq 0$, we define
\begin{align*}
&S_t=\{0\leq i\leq n-1:\eta_t(i)=0\}, ~I_t=\{0\leq i\leq n-1:\eta_t(i)=1\}\\
\text{~and~} &R_t=\{0\leq i\leq n-1:\eta_t(i)=-1\}
\end{align*}
as the sets of vertices in state $0,1,-1$ respectively, then we can identify $\eta_t$ with $(S_t,~I_t,~R_t)$. When the random recovery rates $\{\xi(j)\}_{0\leq j<+\infty}$ and edge weights $\{\rho(i,j)\}_{0\leq i\neq j<+\infty}$ are given, the transition rates of $\{\eta_t\}_{t\geq 0}$ are defined as follows. For any $t\geq 0$ and $0\leq i\leq n-1$,
\begin{equation}\label{equ 1.1 transition rate}
(S_t, ~I_t, ~R_t)\rightarrow
\begin{cases}
(S_t,~I_t\setminus\{i\},~R_t\cup \{i\})\text{~at rate~} \xi(i) & \text{~if~}i\in I_t,\\
(S_t\setminus\{i\}, ~I_t\cup \{i\}, ~R_t)\text{~at rate~}\frac{\lambda}{n}\sum\limits_{j:j\in I_t}\rho(i,j) & \text{~if~}i\in S_t,
\end{cases}
\end{equation}
where $\lambda$ is a positive parameter called the infection rate.

The SIR model describes the spread of an infection on the graph. Vertices in $S_t$ are susceptible to the disease. Vertices in $I_t$ are infective while vertices in $R_t$ are removed. A removed vertex is frozen in this state and will never be infected again. An infective vertex $j$ waits for an exponential time with rate $\xi(j)$ to become removed. A susceptible vertex $i$ is infected by an infective one $j$ at rate proportional to the weight $\rho(i,j)$ on the edge $e(i,j)$.

There are two main types of epidemic models. One is the SIR model which we are concerned with in this paper. The other one is the SIS model, where an infective vertex can become healthy and then be infected again. The stochastic SIS model is also named as the contact process. Readers can see Chapter 6 of \cite{Lig1985} and Part one of \cite{Lig1999} for a survey of the study of the contact process. Recently, the study of the two types of epidemic models in random environments has been a popular topic. We are inspired a lot by related works such as \cite{ Ber2011, Bra1991, Cha2009, Chen2009, Lig1992, Pastor2001a, Pastor2001b, Pet2011, Wang2007, Wang2010, Wang2012, Xue2015, Xue2016, Xue2016b, Yao2012} and so on.

\section{Main result}\label{section two}
In this section we give the main result of this paper. First we introduce some notations and definitions. We assume that the recovery rates $\{\xi(j)\}_{0\leq j<+\infty}$  and edge weights $\{\rho(i,j)\}_{0\leq i\neq j<+\infty}$ are defined under the probability space $\{\Omega,\mathcal{F},\mu\}$. The expectation operator with respect to $\mu$ is denoted by ${\rm E}_\mu$. For any $\omega\in \Omega$, $n\geq 1$ and $\lambda>0$, we denote by $P_{\lambda,n}^\omega$ the probability measure of the process $\{\eta_t\}_{t\geq 0}$ with recovery rates $\{\xi(j,\omega)\}_{0\leq j\leq n-1}$, edge weights $\{\rho(i,j,\omega)\}_{0\leq i\neq j\leq n-1}$ and infection rate $\lambda$ on $C_n$. $P_{\lambda,n}^{\omega}$ is called the quenched measure. We denote by ${\rm E}_{\lambda,n}^\omega$ the expectation operator with respect to $P_{\lambda,n}^\omega$. For each $n\geq 1$, we define
\[
P_{\lambda,n}(\cdot)={\rm E}_\mu\big(P_{\lambda,n}^\omega(\cdot)\big).
\]
$P_{\lambda,n}$ is called the annealed measure. The expectation operator with respect to $P_{\lambda,n}$ is denoted by ${\rm E}_{\lambda,n}$. We assume that at $t=0$, only $0$ is infective and other vertices are susceptible. That is to say,
\[
(S_0,~I_0,~R_0)=(C_n\setminus\{0\},~\{0\},~\emptyset).
\]
According to the transition rates of $\{\eta_t\}_{t\geq 0}$, $R_{t_1}\subseteq R_{t_2}$ for any $t_1<t_2$. Therefore, it is reasonable to define that
\begin{equation}\label{equ 2.1 final removed set}
r_\infty=\lim_{t\rightarrow+\infty}|R_t|,
\end{equation}
where we use $|\cdot|$ to denote the cardinality of a set. According to our assumptions, since $|R_0|=0$,
\begin{equation}\label{equ 2.2}
r_\infty=|\{0\leq i\leq n-1:i\in I_t \text{~for some~}t>0\}|,
\end{equation}
which is the number of vertices that have ever been infective. We write $r_\infty$ as $r_\infty^{(n)}$ when we need to point out that the process is defined on the graph $C_n$.

Now we can give our main result. We obtain that the process $\{\eta_t\}_{t\geq 0}$ performs the following phase transition.
\begin{theorem}\label{theorem 2.1 main}
Let $\lambda_c=\frac{1}{{\rm E}\rho{\rm E}\frac{1}{\xi}}$ and $r_\infty$ be defined as in Equation \eqref{equ 2.1 final removed set}, then for any $\lambda<\lambda_c$ and $\epsilon>0$,
\begin{equation}\label{equ 2.3 subcritical}
\lim_{n\rightarrow+\infty}P_{\lambda,n}\big(\frac{r_\infty^{(n)}}{n}\geq \epsilon\big)=0,
\end{equation}
while for any $\lambda>\lambda_c$, there exist $c(\lambda)>0$ and $b(\lambda)>0$ such that
\begin{equation}\label{equ 2.4 supercritical}
\inf_{n\geq 1} P_{\lambda,n}\Big(\frac{r_\infty^{(n)}}{n}\geq c(\lambda)\Big)~\geq b(\lambda).
\end{equation}
\end{theorem}

According to Theorem \ref{theorem 2.1 main}, when $\lambda<\lambda_c$, the final proportion of vertices that have ever been infective converges to $0$ in distribution while when $\lambda>\lambda_c$, the probabilities that a positive proportion of vertices have ever been infective are uniformly bounded from below by a positive constant.

\proof[Remark 1]

Note that even if $\lambda>\lambda_c$, there is no $c(\lambda)$ which makes $\lim_{n\rightarrow+\infty}P_{\lambda,n}(\frac{r_\infty^{(n)}}{n}\geq c(\lambda))=1$. As a result, Theorem \ref{theorem 2.1 main} can not be strengthened to a zero-one law. The reason of the nonexistence of such $c(\lambda)$ is that
\[
P_{\lambda,n}(r_\infty=1)=P_{\lambda,n}(0\text{~does not infect any other~})=P_{\lambda,n}(U\leq \inf_{1\leq i\leq n-1}T_i),
\]
where $U$ is an exponential time with rate $\xi(0)$ while $T_i$ is an exponential time with rate $\frac{\lambda}{n}\rho(0,i)$ for each $i$ and all these exponential times are independent under $P_{\lambda,n}^{\omega}$. Therefore,
\[
P_{\lambda,n}(r_\infty=1)={\rm E_\mu}\Big(\frac{\xi(0)}{\xi(0)+\frac{\lambda}{n}\sum\limits_{i=1}^{n-1}\rho(0,i)}\Big)\rightarrow {\rm E}\frac{\xi}{\xi+\lambda{\rm E}\rho}
\]
according to the law of large numbers. Then,
\[
\limsup_{n\rightarrow+\infty}P_{\lambda,n}\Big(\frac{r_\infty^{(n)}}{n}\geq c(\lambda)\Big)\leq  \lim_{n\rightarrow+\infty}P_{\lambda,n}\big(r_\infty^{(n)}> 1\big)=\lambda {\rm E}\rho{\rm E}\big(\frac{1}{\xi+\lambda{\rm E}\rho}\big)<1
\]
for any $c(\lambda)>0$.

\qed

\proof[Remark 2]
The intuitive explanation of $\lambda_c=\frac{1}{{\rm E}\rho{\rm E}\frac{1}{\xi}}$ is as follows. The probability that an infective vertex infects a given susceptible one before becoming removed is ${\rm E}\big(\frac{\frac{\lambda}{n}\rho}{\frac{\lambda}{n}\rho+\xi}\big)\approx \frac{1}{n}\frac{\lambda}{\lambda_c}$ for large $n$. Then, $r_\infty$ with order $O(n)$ is somewhat similar with that the Erd\H{o}s-R\'{e}nyi graph $G(n,\frac{1}{n}\frac{\lambda}{\lambda_c})$ has the largest component with order $O(n)$. Erd\H{o}s-R\'{e}nyi graph $G(n,\frac{\mu}{n})$ has critical value $\mu=1$ (see Section 3 of \cite{Hofstad2013}), hence phase transition should occur when $\lambda=\lambda_c$ for our model according to the above non-rigorous comparison.

\qed

\proof[Remark 3]
When $\xi=\rho\equiv 1$, our model reduces to the classic SIR model on complete graphs, where vertices in the same state are equal. Then according to the classic theory of density dependent population process (see Chapter 11 of \cite{Ethier1986}), $(|S_t|/n,~|I_t|/n, ~|R_t|/n)$ converges weakly to the solution $(\widetilde{s}_t,\widetilde{i}_t,\widetilde{r}_t)$ of the following  deterministic ODE.
\[
\begin{cases}
&\frac{d}{dt}\widetilde{s}_t=-\lambda\widetilde{i}_t\widetilde{s}_t,\\
&\frac{d}{dt}\widetilde{i}_t=\widetilde{i}_t(\lambda \widetilde{s}_t-1),\\
&\frac{d}{dt}\widetilde{r}_t=\widetilde{i}_t,
\end{cases}
\]
where $\lambda_c=1$ is the maximum of $\lambda$ that makes $\frac{d}{dt}\widetilde{i}_t<0$ for any initial condition and $t>0$, which implies that the infection can not spread in large scale when $\lambda<\lambda_c$.

\qed

\proof[Remark 4]
In \cite{Pet2011}, Peterson considers the SIS model with random transition rates (in detail, random vertex weights) on complete graphs. It is shown in \cite{Pet2011} that the SIS model there has the critical value which equals the inverse of the second moment of the vertex weight. The phases of the SIS model are divided according to time when the process dies out. In the supercritical case, the process survives at time $\exp\{O(n)\}$ with high probability while in the subcritical case, the process dies out at time $O(\log n)$ with high probability. However, the spread of the epidemic with the SIR type will not last for a long time since there is no repeated infections. As a result, we divide the phases of the SIR model according to the proportion of vertices that have ever been infected.

\qed

The proof of Theorem \ref{theorem 2.1 main} is divided into two sections. In Section \ref{section three}, we give the proof of Equation \eqref{equ 2.3 subcritical}. The core idea is to show that $\{{\rm E}_{\lambda,n}r_\infty^{(n)}\}_{n\geq 1}$ are uniformly bounded from above when $\lambda<\lambda_c$. To do so, we will show that the mean of the number of paths with length $k$ through which the infection spreads is about $(\lambda/\lambda_c)^k$ for large $k$. In
Section \ref{section four}, we give the proof of Equation \eqref{equ 2.4 supercritical}. The core idea is to show that ${\rm E}_{\lambda,n}r_\infty^{(n)}=O(n)$ for large $n$ when $\lambda>\lambda_c$, which is equivalent to that $\{P_{\lambda,n}\big(1\in I_t\text{~for some~}t>0\big)\}_{n\geq 1}$ are uniformly bounded from below. We will consider the paths with length $O(\log n)$ through which the infection spreads from $0$ to $1$. H\"{o}lder inequality will be crucial for us to bound from below the probability of the existence of such paths.

\section{Proof of Equation \eqref{equ 2.3 subcritical}}\label{section three}

In this section we give the proof of Equation \eqref{equ 2.3 subcritical}. First we introduce some notations and definitions. For any $0\leq i,j\leq n-1$, let $T_i$ be an exponential time with rate $\xi(i)$ and $U(i,j)$ be an exponential time with rate $\frac{\lambda}{n}\rho(i,j)$. Note that here we care about the order of $i$ and $j$, hence $U(i,j)\neq U(j,i)$. We assume that all these exponential times are defined under the quenched measure $P_{\lambda,n}^\omega$ and are independent under $P_{\lambda,n}^{\omega}$. Note that $U(i,j)$ and $U(j,i)$ are positive correlated under the annealed measure $P_{\lambda,n}$ since they are both with rate $\frac{\lambda}{n}\rho(i,j)$. Intuitively, $T(i)$ is the time interval which $i$ waits for to become removed after it is infected while $U(i,j)$ is the time interval which $i$ waits for to infect $j$ after $i$ is infected. For $0\leq k\leq n-1$, we denote by $L_k$ the set of self-avoiding paths on $C_n$ starting at $0$ with length $k$. That is to say,
\[
L_k=\big\{\vec{l}=(l_0,l_1,\ldots,l_k)\in \{0,1,\ldots,n-1\}^k:l_0=0, l_i\neq l_j\text{~for any~}0\leq i<j\leq k\big\}.
\]
Now we give the proof of Equation \eqref{equ 2.3 subcritical}.

\proof

According to the transition rates of our process, $i\in I_t$ for some $t>0$ when and only when there exist $k>0$ and $\vec{l}=(l_0,l_1,\ldots,l_k)\in L_k$ with $l_k=i$ such that $U(l_j,l_{j+1})\leq T(l_j)$ for $0\leq j\leq k-1$ in the sense of coupling. For $\vec{l}=(l_0,l_1,\ldots,l_k)\in L_k$, we denote by $A_{\vec{l}}$ the event that $U(l_j,l_{j+1})\leq T(l_j)$ for $0\leq j\leq k-1$. Then, according to the above analysis,
\begin{equation}\label{equ 3.1}
P_{\lambda,n}(i\in I_t,\text{~for some~}t>0)\leq \sum_{k=0}^{n-1}\sum_{\vec{l}\in L_k,\atop l_k=i}P_{\lambda,n}(A_{\vec{l}}),
\end{equation}
where $\vec{l}=(l_0,l_1,\ldots,l_k)$. By Equations \eqref{equ 2.2} and \eqref{equ 3.1},
\begin{align}\label{equ 3.2}
{\rm E}_{\lambda,n}R_{\infty}^{(n)}&=\sum_{i\in C_n}P_{\lambda,n}(i\in I_t,\text{~for some~}t>0)\notag\\
&=\sum_{k=0}^{n-1}\sum_{\vec{l}\in L_k}P_{\lambda,n}(A_{\vec{l}}).
\end{align}

Since $\vec{l}$ is self-avoiding while the recovery rates and edge weights are independent under the annealed measure,
\begin{align}\label{equ 3.3}
P_{\lambda,n}(A_{\vec{l}})&={\rm E}_\mu\Big(\prod_{j=0}^{k-1}P_{\lambda,n}^{\omega}\big(U(l_j,l_{j+1})\leq T(j)\big)\Big)\notag\\
&={\rm E}_\mu\Big(\prod_{j=0}^{k-1}\frac{\frac{\lambda}{n}\rho(l_j,l_{j+1})}{\frac{\lambda}{n}\rho(l_j,l_{j+1})+\xi(l_j)}\Big)=\big({\rm E}\frac{\frac{\lambda}{n}\rho}{\frac{\lambda}{n}\rho+\xi}\big)^k\leq \frac{\lambda^k}{n^k\lambda_c^k}.
\end{align}

By Equations \eqref{equ 3.2} and \eqref{equ 3.3}, for $\lambda<\lambda_c$,
\begin{equation}\label{equ 3.4}
{\rm E}_{\lambda,n}R_{\infty}^{(n)}\leq \sum_{k=0}^{+\infty}\big(\frac{\lambda}{\lambda_c}\big)^k=\frac{\lambda_c}{\lambda_c-\lambda},
\end{equation}
since $|L_k|=n(n-1)\ldots(n-k)\leq n^k$.

By Equation \eqref{equ 3.4} and Chebyshev's inequality, for $\lambda<\lambda_c$ and $\epsilon>0$,
\[
P_{\lambda,n}\big(\frac{R_\infty^{(n)}}{n}\geq \epsilon\big)\leq \frac{\lambda_c}{\epsilon n(\lambda_c-\lambda)}\rightarrow 0
\]
as $n\rightarrow+\infty$.

\qed

\section{Proof of Equation \eqref{equ 2.4 supercritical}}\label{section four}

In this section we give the proof of Equation \eqref{equ 2.4 supercritical}. First we give some notations and definitions. For $\lambda>\lambda_c$, we fix $\delta=\delta(\lambda)>0$ and $\theta=\theta(\lambda)>0$ such that
\[
\frac{\lambda(1-\delta)}{\lambda_c}>1\text{~and~} \big(\frac{\lambda(1-\delta)}{\lambda_c}\big)^\theta>e.
\]
For sufficiently large $n$ such that $\theta \log n<n$, we denote by $\{S_k^{(n)}\}_{0\leq k\leq \lfloor\theta \log n\rfloor}$ the self-avoiding random walk on $C_n$ with $S_0^{(n)}=0$ and $S_{\lfloor\theta\log n\rfloor}^{(n)}=1$. That is to say, for each $0\leq k\leq \lfloor\theta \log n\rfloor-1$ and any $\vec{l}=(l_0,l_1,\ldots,l_k)\in L_k$ with $l_j\neq 1$ for $1\leq j\leq k$,
\[
P(S^{(n)}_j=l_j\text{~for all~}0\leq j\leq k, ~S^{(n)}_{\lfloor\theta \log n\rfloor}=1)=\frac{1}{(n-2)(n-3)\ldots(n-k-1)}.
\]
We denote by $\{\widetilde{S}^{(n)}_k\}_{0\leq k\leq \lfloor \theta \log n\rfloor}$ an independent copy of $\{S_k^{(n)}\}_{0\leq k\leq \lfloor \theta \log n\rfloor}$. We assume that $\{S_k^{(n)}\}_{0\leq k\leq \lfloor \theta \log n\rfloor}$ and $\{\widetilde{S}^{(n)}_k\}_{0\leq k\leq \lfloor \theta \log n\rfloor}$ are defined under the probability space $(\widetilde{\Omega},\widetilde{\mathcal{F}},\widetilde{P})$. The expectation operator with respect to $\widetilde{P}$ is denoted by $\widetilde{{\rm E}}$. When there is no misunderstanding, we write $S_k^{(n)}$ and $\widetilde{S}_k^{(n)}$ as $S_k$ and $\widetilde{S}_k$.

For later use, we introduce the following definitions. For $\vec{x}=(x_0,x_1,\ldots,x_{\lfloor\theta \log n\rfloor})\in L_{\lfloor\theta \log n\rfloor}$ and
$\vec{y}=(y_0,y_1,\ldots,y_{\lfloor\theta \log n\rfloor})\in L_{\lfloor\theta \log n\rfloor}$, we define
\[
D(\vec{x},\vec{y})=\{0\leq i\leq \lfloor\theta \log n\rfloor-1: y_i=x_j\text{~for some~}0\leq j\leq \lfloor\theta \log n\rfloor-1\}
\]
and
\[
F(\vec{x},\vec{y})=\{0\leq i\leq \lfloor\theta \log n\rfloor-1: y_i=x_j \text{~and~}y_{i+1}=x_{j+1}\text{~for some~}0\leq j\leq \lfloor\theta \log n\rfloor-1\}.
\]
Note that $|D(\vec{x},\vec{y})|$ is the number of vertices both $\vec{x}$ and $\vec{y}$ visit while $|F(\vec{x},\vec{y})|$ is the number of edges $\vec{y}$ goes through in the direction as that of $\vec{x}$. Since $\vec{x}$ and $\vec{y}$ are self-avoiding, if $i\in D(\vec{x},\vec{y})$, then there exists an unique $j$ such that $y_i=x_j$.

We use $S$ and $\widetilde{S}$ to denote the random paths $(S_0,\ldots, S_{\lfloor\theta \log n\rfloor})$ and $(\widetilde{S}_0,\ldots,\widetilde{S}_{\lfloor\theta\log n\rfloor})$, which belong to $L_{\lfloor\theta\log n\rfloor}$, then we have the following lemma which is crucial for us to prove Equation \eqref{equ 2.4 supercritical}.
\begin{lemma}\label{lemma 4.1}
For any $\lambda>\lambda_c$ and sufficiently large $n$,
\[
P_{\lambda,n}(1\in I_t\text{~for some~}t>0)\geq \frac{1}{\widetilde{\rm E}\Big(\big(\frac{n\lambda_c}{(1-\delta)\lambda}\big)^{|F(S,\widetilde{S})|}\big(\frac{2{\rm E}\frac{1}{\xi^2}}{({\rm E}\frac{1}{\xi})^2(1-\delta)^2{\rm E}\rho}\big)^{|D(S,\widetilde{S})\setminus F(S,\widetilde{S})|}\Big)}.
\]
\end{lemma}

\proof

Let
\[
\widehat{L}=\big\{\vec{l}=(l_0,l_1,\ldots,l_{\lfloor \theta\log n\rfloor})\in L_{\lfloor \theta\log n\rfloor}:l_{\lfloor \theta\log n\rfloor}=1\big\}
\]
be the set of all the possible paths of $S$ and $\widetilde{S}$, then for each $\vec{l}\in \widehat{L}$, we denote by $\chi_{\vec{l}}$ the indicator function of the event $A_{\vec{l}}$, where $A_{\vec{l}}$ is defined as in Section \ref{section three}. According to the transition rates of the process,
\begin{equation}\label{equ 4.1}
\{1\in I_t\text{~for some~}t>0\}\supseteq \{\sum_{\vec{l}\in \widehat{L}}\chi_{\vec{l}}>0\}.
\end{equation}
By Equation \eqref{equ 4.1} and H\"{o}lder's inequality,
\begin{equation}\label{equ 4.2}
P_{\lambda,n}\big(1\in I_t\text{~for some~}t>0\big)\geq \frac{\big({\rm E}_{\lambda,n}\sum\limits_{\vec{l}\in \widehat{L}}\chi_{\vec{l}}\big)^2}{{\rm E}_{\lambda,n}\Big(\big(\sum\limits_{\vec{l}\in \widehat{L}}\chi_{\vec{l}}\big)^2\Big)}
=\frac{\sum\limits_{\vec{x}\in\widehat{L}}\sum\limits_{\vec{y}\in \widehat{L}}P_{\lambda,n}(A_{\vec{x}})P_{\lambda,n}(A_{\vec{y}})}
{\sum\limits_{\vec{x}\in\widehat{L}}\sum\limits_{\vec{y}\in \widehat{L}}P_{\lambda,n}(A_{\vec{x}}\cap A_{\vec{y}})}.
\end{equation}

According to the spatial homogeneity of our model under the annealed measure,
\[P_{\lambda,n}(A_{\vec{x}})=P_{\lambda,n}(A_{\vec{y}})\]
for any $\vec{x},\vec{y}\in \widehat{L}$. Then by Equation \eqref{equ 4.2},
\begin{align}\label{equ 4.3}
P_{\lambda,n}\big(1\in I_t\text{~for some~}t>0\big)&\geq \frac{1}{\frac{1}{|\widehat{L}|^2}\sum\limits_{\vec{x}\in\widehat{L}}\sum\limits_{\vec{y}\in \widehat{L}}\frac{P_{\lambda,n}(A_{\vec{x}}\cap A_{\vec{y}})}{P_{\lambda,n}(A_{\vec{x}})P_{\lambda,n}(A_{\vec{y}})}}\notag\\
&=\frac{1}{\frac{1}{|\widehat{L}|^2}\sum\limits_{\vec{x}\in\widehat{L}}\sum\limits_{\vec{y}\in \widehat{L}}\frac{P_{\lambda,n}(A_{\vec{y}}|A_{\vec{x}})}{P_{\lambda,n}(A_{\vec{y}})}}.
\end{align}

Now we deal with $\frac{P_{\lambda,n}(A_{\vec{y}}|A_{\vec{x}})}{P_{\lambda,n}(A_{\vec{y}})}$. Let $\{T(i)\}_{0\leq i\leq n-1}$ and $\{U(i,j)\}_{0\leq i,j\leq n-1}$ be defined as in Section \ref{section three}, then $T(i)$ is an exponential time with rate $\xi(i)$ for each $i$ while $U(i,j)$ is an exponential time with rate $\frac{\lambda}{n}\rho(i,j)$ for any $i\neq j$.  Note that $\vec{x}$ and $\vec{y}$ are self-avoiding while the recovery rates and edge weights are independent, hence for each $i\not \in D(\vec{x},\vec{y})$, both the numerator and the denominator of $\frac{P_{\lambda,n}(A_{\vec{y}}|A_{\vec{x}})}{P_{\lambda,n}(A_{\vec{y}})}$ have the factor  $P_{\lambda,n}\big(U(y_i, y_{i+1})\leq T(y_i)\big)$, which can be cancelled. For $i\in F(\vec{x},\vec{y})$, the denominator of $\frac{P_{\lambda,n}(A_{\vec{y}}|A_{\vec{x}})}{P_{\lambda,n}(A_{\vec{y}})}$ has the factor $P_{\lambda,n}\big(U(y_i,y_{i+1})\leq T(y_i)\big)$ but the numerator does not. This is because $(y_i,y_{i+1})=(x_j,x_{j+1})$ for some $j$, then $U(y_i,y_{i+1})\leq T(y_i)$ occurs with probability one conditioned on $A_{\vec{x}}$. Therefore, for each $i\in F(\vec{x},\vec{y})$, there is a factor
\[
\frac{1}{P_{\lambda,n}\big(U(y_i, y_{i+1})\leq T(y_i)\big)}=\frac{1}{\frac{\lambda}{n}{\rm E}\big(\frac{\rho}{\frac{\lambda}{n}\rho+\xi}\big)}
\]
in $\frac{P_{\lambda,n}(A_{\vec{y}}|A_{\vec{x}})}{P_{\lambda,n}(A_{\vec{y}})}$. For each $i\in D(\vec{x},\vec{y})\setminus F(\vec{x},\vec{y})$, there exist $m,q,r$ which are different with each other such that $x_j=y_i=m$, $x_{j+1}=q$ and $y_{i+1}=r$ for some $j$. Then, the denominator of $\frac{P_{\lambda,n}(A_{\vec{y}}|A_{\vec{x}})}{P_{\lambda,n}(A_{\vec{y}})}$ has factor $P_{\lambda,n}\big(U(m, r)\leq T(m)\big)$ while the numerator has a factor no more than
$P_{\lambda,n}\big(\widetilde{U}(m, r)\leq T(m)|U(m,q)\leq T(m)\big)$, where $\widetilde{U}(m,r)$ is an exponential time with rate $\frac{\lambda}{n}$ and independent with other exponential times, since $\rho(m,r)\leq 1$ according to our assumption. Note that we replace $U(m,r)$ by $\widetilde{U}(m,r)$ to ignore the correlation between $U(m,r)$ and $U(r,m)$ under the annealed measure when $x_{j-1}=r$. Hence $\frac{P_{\lambda,n}(A_{\vec{y}}|A_{\vec{x}})}{P_{\lambda,n}(A_{\vec{y}})}$ has a factor no more than
\[
\frac{P_{\lambda,n}\big(\widetilde{U}(m, r)\leq T(m)|U(m,q)\leq T(m)\big)}{P_{\lambda,n}\big(U(m, r)\leq T(m)\big)},
\]
which is no more than
\[
\frac{2{\rm E}\frac{1}{\xi^2}{\rm E}\rho}{\big({\rm E}\frac{\rho}{\xi+\frac{\lambda}{n}\rho}\big)^2}
\]
by direct calculation. According to the above analysis,
\begin{equation}\label{equ 4.4}
\frac{P_{\lambda,n}(A_{\vec{y}}|A_{\vec{x}})}{P_{\lambda,n}(A_{\vec{y}})}\leq \Big(\frac{1}{\frac{\lambda}{n}{\rm E}\big(\frac{\rho}{\frac{\lambda}{n}\rho+\xi}\big)}\Big)^{|F(\vec{x},\vec{y})|}\Big(\frac{2{\rm E}\frac{1}{\xi^2}{\rm E}\rho}{\big({\rm E}\frac{\rho}{\xi+\frac{\lambda}{n}\rho}\big)^2}\Big)^{|D(\vec{x},\vec{y})\setminus F(\vec{x},\vec{y})|}.
\end{equation}
According to Dominated Convergence Theorem, for sufficiently large $n$,
\begin{equation}\label{equ 4.5}
{\rm E}\big(\frac{\rho}{\frac{\lambda}{n}\rho+\xi}\big)\geq \frac{1-\delta}{\lambda_c}.
\end{equation}
Then, Lemma \ref{lemma 4.1} follows from Equations \eqref{equ 4.3}, \eqref{equ 4.4} and \eqref{equ 4.5}, since
\[
\widetilde{\rm E}\Big(A^{|F(S,\widetilde{S})|}~B^{|D(S,\widetilde{S})\setminus F(S,\widetilde{S})|}\Big)
=\frac{1}{|\widehat{L}|^2}\sum\limits_{\vec{x}\in\widehat{L}}\sum\limits_{\vec{y}\in \widehat{L}}A^{|F(\vec{x},\vec{y})|}~B^{|D(\vec{x},\vec{y})\setminus F(\vec{x},\vec{y})|}
\]
for any $A,B>0$.

\qed

According to Lemma \ref{lemma 4.1}, we can show that $\{P_{\lambda,n}\big(1\in I_t\text{~for some~}t>0\big)\}_{n\geq 1}$ are uniformly bounded by a positive constant from below.
\begin{lemma}\label{lemma 4.2}
For each $\lambda>\lambda_c$, there exists $d(\lambda)>0$ such that
\[
\inf_{n\geq 1} P_{\lambda,n}\big(1\in I_t\text{~for some~}t>0\big)~\geq d(\lambda).
\]
\end{lemma}

We give the proof of Lemma \ref{lemma 4.2} at the end of this section. First we show how to use Lemma \ref{lemma 4.2} to prove
Equation \eqref{equ 2.4 supercritical}.

\proof[Proof of Equation \eqref{equ 2.4 supercritical}]

According to Equation \eqref{equ 2.2}, Lemma \ref{lemma 4.2} and the spatial homogeneity of the model under the annealed measure $P_{\lambda,n}$, for any $\lambda>\lambda_c$,
\begin{equation}\label{equ 4.6}
{\rm E}_{\lambda,n}r_\infty^{(n)}=1+(n-1)P_{\lambda,n}(1\in I_t\text{~for some~}t>0)\geq 1+(n-1)d(\lambda)\geq nd(\lambda)
\end{equation}
for each $n\geq1$, where $d(\lambda)$ is defined as in Lemma \ref{lemma 4.2}. Note that $1$ on the right-hand side of the above equality is the probability that $0$ will be removed. By Equation \eqref{equ 4.6},
\begin{align*}
nd(\lambda)&\leq {\rm E}_{\lambda,n}r_\infty^{(n)} \leq \frac{nd(\lambda)}{2}P_{\lambda,n}\Big(r_\infty^{(n)}
\leq \frac{d(\lambda)n}{2}\Big)+nP_{\lambda,n}\Big(r_\infty^{(n)}\geq \frac{d(\lambda)n}{2}\Big)\\
&\leq \frac{nd(\lambda)}{2}+nP_{\lambda,n}\Big(r_\infty^{(n)}\geq \frac{d(\lambda)n}{2}\Big)
\end{align*}
and hence
\[
P_{\lambda,n}\Big(\frac{r_\infty^{(n)}}{n}\geq \frac{d(\lambda)}{2}\Big)\geq \frac{d(\lambda)}{2}>0
\]
for each $n\geq 1$. Let $c(\lambda)=b(\lambda)=\frac{d(\lambda)}{2}$, then the proof is complete.

\qed

At last we give the proof of Lemma \ref{lemma 4.2}. Let us explain the intuitive idea of proof first. For random walk $S$ and $\widetilde{S}$, each step has about $n$ choices (in detail, at least $n-\lfloor\theta\log n\rfloor$ choices), hence $P(S_{j+1}=\widetilde{S}_{i+1}|S_j=\widetilde{S}_i)\approx \frac{1}{n}$. There are $O(\log n)$ vertices on the paths $S$ and $\widetilde{S}$, hence $P(\widetilde{S}_{i}=S_j\text{~for some~}i\text{~and~}j|S)\leq  O(\frac{(\log n)^2}{n})$. As a result, according to the strong Markov property, we expect that
\[
P\big(|F(S,\widetilde{S})|=k,~|D(S,\widetilde{S})\setminus F(S,\widetilde{S})|=l\big)\leq (\frac{1}{n})^k\big(\frac{C_0(\log n)^2}{n}\big)^l
\]
for $k,l>0$, where $C_0$ is a constant.

Then, the denominator of the right-hand side of the inequality in Lemma \ref{lemma 4.1} is bounded from above by
\begin{equation}\label{equ 4.7}
\sum_{l\geq0,k\geq 0}(C_1 n)^k(C_2)^l(\frac{1}{n})^k\big(\frac{C_0(\log n)^2}{n}\big)^l<+\infty
\end{equation}
for large $n$, where $C_1=\frac{\lambda_c}{\lambda(1-\delta)}$ and $C_2>0$ are constants. Note that $C_1<1$, which only holds in the case where $\lambda>\lambda_c$, is crucial for the above upper bound to be finite. Readers may wonder why we are concerned with the infection paths with length $\lfloor\theta\log n\rfloor$ instead of shorter paths to execute the above analysis. The reason relies on the following detail. Since $S$ and $\widetilde{S}$ both start at $0$ and end at $1$, when $\widehat{k}$ is the length of $S$,
\[
P\big(|F(S,\widetilde{S})|=\widehat{k}\big)\approx\big(\frac{1}{n}\big)^{\widehat{k}-1},
\]
not $(\frac{1}{n})^{\widehat{k}}$, which is not consistent with the above analysis and is the only exception. This exception generates a term $nC_1^{\widehat{k}}$ in Equation \eqref{equ 4.7}. To make this term converge to $0$ as $n\rightarrow+\infty$, we need the length $\widehat{k}$ with order $O(\log n)$.

\proof[Proof of Lemma \ref{lemma 4.2}]

In this proof, we use $C_1$ to denote $\frac{\lambda_c}{(1-\delta)\lambda}$ and $C_2$ to denote $\frac{2{\rm E}\frac{1}{\xi^2}}{({\rm E}\frac{1}{\xi})^2(1-\delta)^2{\rm E}\rho}$. By Lemma \ref{lemma 4.1}, we only need to show that
\begin{equation}\label{equ 4.8}
\sup_{n\geq 1}\widetilde{\rm E}\Big(\big(nC_1\big)^{|F(S,\widetilde{S})|}~C_2^{|D(S,\widetilde{S})\setminus F(S,\widetilde{S})|}\Big)~<+\infty.
\end{equation}
Obviously, we only need to deal with sufficiently large $n$, so we assume that
\[
n-\lfloor\theta \log n\rfloor\geq (1-\delta_2)n,
\]
where $\delta_2>0$ is sufficiently small such that
\[
\frac{C_1}{1-\delta_2}<1\text{~and~}\big(\frac{1-\delta_2}{C_1}\big)^\theta>e.
\]
Note that $C_1<1$ and $\big(\frac{1}{C_1}\big)^\theta>e$ as we have introduced at the beginning of this section, so there exists such $\delta_2$. We define
\[
\tau=\inf\{0<i<\lfloor \theta\log n\rfloor:~S_i\neq\widetilde{S}_i\}
\]
as the first moment when $S$ does not collide with $\widetilde{S}$ and define $\inf \emptyset=+\infty$, then
\begin{equation}\label{equ 4.9}
\big(nC_1\big)^{|F(S,\widetilde{S})|}~C_2^{|D(S,\widetilde{S})\setminus F(S,\widetilde{S})|}=\big(nC_1\big)^{\lfloor\theta\log n\rfloor}.
\end{equation}
on the event $\{\tau=+\infty\}$. According to definition of $S$,
\begin{equation}\label{equ 4.10}
\widetilde{P}\big(\tau=+\infty\big)\leq \big(\frac{1}{n-\lfloor\theta \log n\rfloor}\big)^{\lfloor\theta \log n\rfloor-1}
\leq \big(\frac{1}{n(1-\delta_2)}\big)^{\lfloor\theta \log n\rfloor-1}.
\end{equation}

On the event $\{\tau<+\infty\}$, we introduce the following notations,
\[
\sigma_1=\tau-1=\sup\{i:S_j=\widetilde{S}_j\text{~for all~}0\leq j\leq i\},
\]
and
\[
\sigma_2=\lfloor\theta \log n\rfloor-\inf\{i:S_j=\widetilde{S}_j\text{~for all~}i\leq j\leq \lfloor\theta \log n\rfloor\},
\]
then $\sigma_1<\lfloor\theta \log n\rfloor-\sigma_2$ on $\{\tau<+\infty\}$. On the event $\{\tau<+\infty\}$, we define
\[
\kappa=|\{\sigma_1<i<\lfloor\theta \log n\rfloor-\sigma_2:i\in D(S,\widetilde{S})\}|.
\]
For $0\leq l\leq \kappa$, we define $t(0)=\sigma_1$ and
\[
t(l)=\inf\{i>t(l-1):i\in D(S,\widetilde{S})\}
\]
for $1\leq l\leq \kappa$. For $0\leq l\leq \kappa$, we define
\[
H(l)=
\begin{cases}
nC_1 &\text{~if\quad} t(l)\in F(S,\widetilde{S}),\\
C_2 &\text{~if\quad} t(l)\in D(S,\widetilde{S})\setminus F(S,\widetilde{S}),
\end{cases}
\]
then on the event $\tau<+\infty$,
\begin{equation}\label{equ 4.11}
\big(nC_1\big)^{|F(S,\widetilde{S})|}~C_2^{|D(S,\widetilde{S})\setminus F(S,\widetilde{S})|}=\big(nC_1\big)^{(\sigma_1+\sigma_2)}\prod_{l=0}^{\kappa}H(l),
\end{equation}
since $\{0,1,\ldots,\sigma_1-1\}\bigcup\{\lfloor\theta\log n\rfloor-\sigma_2,\lfloor\theta\log n\rfloor-\sigma_2+1,\ldots,\lfloor\theta\log n\rfloor-1\}~\subseteq F(S,\widetilde{S})$. Note that according to our definitions, $H(0)=H(\kappa)=C_2$. By the definition of $S$ and $\widetilde{S}$, for any integers $l,k\geq 0$,
\begin{equation}\label{equ 4.12}
\widetilde{P}(\sigma_1=l,\sigma_2=k)\leq \big(\frac{1}{n-\lfloor\theta \log n\rfloor}\big)^{l+k}\leq \big(\frac{1}{n(1-\delta_2)}\big)^{l+k}.
\end{equation}
For $l\geq1$, on the event $\{\kappa\geq l\}$, $H(l)=C_1n$ when and only when $\widetilde{S}_{t(l)+1}=S_{j+1}$, where $j$ is the unique vertex such that $\widetilde{S}_{t(l)}=S_j$. Therefore, conditioned on $S$, $\{t(j)\}_{0\leq j\leq l}$ and $\{\widetilde{S}_j\}_{0\leq j\leq t(l)}$, the probability that $H(l)=C_1n$ is bounded from above by
\[
\max_{1\leq j\leq \lfloor\theta \log n\rfloor-1,\atop k\neq 0,1}\widetilde{P}(\widetilde{S}_j=k)\leq \frac{1}{n-\lfloor\theta \log n\rfloor}\leq \frac{1}{n(1-\delta_2)}.
\]
Note that $H(l)$ can not be $nC_1$ when $t(l)=\lfloor\theta \log n\rfloor-\sigma_2-1$ according to the definition of $\sigma_2$. Hence the information of $\sigma_2$ can not enlarge the probability that $H(l)=nC_1$, which fact we use to obtain the above upper bound.

As a result, for $l\geq 1$,
\begin{align}\label{equ 4.13}
&\widetilde{E}\Big(H(l);~H(l)=C_1n,\kappa\geq l+1\Big|\sigma_1,\sigma_2,\{t(i)\}_{0\leq i\leq l},~S,~\{\widetilde{S}_i\}_{0\leq i\leq t(l)}\Big)\\
&\leq \frac{C_1n}{n(1-\delta_2)}=\frac{C_1}{1-\delta_2} \notag
\end{align}
on the event $\{\kappa\geq l\}$. Note that $H(l)=C_1n$ implies that $t(l+1)=t(l)+1$ and as a result $\kappa\geq l+1$. For $l\geq 0$, on the event $\{\kappa\geq l\}$, $\kappa\geq l+1$ when there exists $i>t(l)$ and $j$ such that $\widetilde{S}_i=S_j$. Therefore, conditioned on $S$, $\{t(j)\}_{0\leq j\leq l}$ and $\{\widetilde{S}_j\}_{0\leq j\leq t(l)}$, the probability that $\kappa\geq l+1$ is bounded from above by
\begin{align*}
\sum_{1\leq i,j\leq \lfloor\theta \log n\rfloor}\widetilde{P}(\widetilde{S}_i=S_j|S)&\leq \sum_{1\leq i,j\leq \lfloor\theta \log n\rfloor}\max_{k\neq 0,1}\widetilde{P}(\widetilde{S}_i=k)\\
&\leq \frac{\lfloor\theta \log n\rfloor^2}{n-\lfloor\theta \log n\rfloor}\leq \frac{M_2(\log n)^2}{n},
\end{align*}
where $M_2$ is a positive constant which does not depend on $n$. As a result, for each $l\geq 0$,
\begin{equation}\label{equ 4.14}
\widetilde{E}\Big(H(l);~H(l)=C_2,\kappa\geq l+1\Big|\sigma_1,\sigma_2,\{t(i)\}_{0\leq i\leq l},~S,~\{\widetilde{S}_i\}_{0\leq i\leq t(l)}\Big)
\leq \frac{C_2M_2(\log n)^2}{n}
\end{equation}
and
\begin{equation}\label{equ 4.15}
\widetilde{E}\Big(H(l);~H(l)=C_2,\kappa=l\Big|\sigma_1,\sigma_2,\{t(i)\}_{0\leq i\leq l},~S,~\{\widetilde{S}_i\}_{0\leq i\leq t(l)}\Big)
\leq C_2
\end{equation}
on the event $\{\kappa\geq l\}$. By Equations \eqref{equ 4.13}, \eqref{equ 4.14}, \eqref{equ 4.15} and strong Markov property,
\begin{equation}\label{equ 4.16}
\widetilde{\rm E}\Big(\prod_{l=0}^{\kappa}H(l)\Big|\sigma_1,\sigma_2\Big)\leq\Big(\sum_{l=0}^{+\infty}\Gamma^l(2,2)\Big)C_2\leq C_2\sum_{l=0}^{+\infty}\big(\frac{C_1}{1-\delta_2}+\frac{C_2M_2(\log n)^2}{n}\big)^l
\end{equation}
since $H(0)=H(\kappa)=C_2$, where $\Gamma$ is a $2\times 2$ matrix such that
\[
\Gamma=
\begin{pmatrix}
\frac{C_1}{1-\delta_2} & \frac{C_2M_2(\log n)^2}{n}\\
\frac{C_1}{1-\delta_2} & \frac{C_2M_2(\log n)^2}{n}
\end{pmatrix}.
\]
Fix $\gamma\in (\frac{C_1}{1-\delta_2},1)$, then
\[
\frac{C_1}{1-\delta_2}+\frac{C_2M_2(\log n)^2}{n}<\gamma
\]
for sufficiently large $n$. Hence by Equation \eqref{equ 4.16}, for sufficiently large $n$,
\begin{equation}\label{equ 4.17}
\widetilde{\rm E}\Big(\prod_{l=0}^{\kappa}H(l)\Big|\sigma_1,\sigma_2\Big)\leq \frac{C_2}{1-\gamma}.
\end{equation}
By Equations \eqref{equ 4.11}, \eqref{equ 4.12} and \eqref{equ 4.17},
\begin{align}\label{equ 4.18}
\widetilde{\rm E}\Big(\big(nC_1\big)^{|F(S,\widetilde{S})|}~C_2^{|D(S,\widetilde{S})\setminus F(S,\widetilde{S})|};~\tau<+\infty\Big)
&\leq \frac{C_2}{1-\gamma}\sum_{l,k\geq 0}(nC_1)^{l+k}\big(\frac{1}{n(1-\delta_2)}\big)^{l+k} \notag\\
&=\frac{C_2}{1-\gamma}\big(\frac{1}{1-\frac{C_1}{1-\delta_2}}\big)^2<+\infty.
\end{align}
By Equations \eqref{equ 4.9} and \eqref{equ 4.10},
\begin{align}\label{equ 4.19}
\widetilde{\rm E}\Big(\big(nC_1\big)^{|F(S,\widetilde{S})|}~C_2^{|D(S,\widetilde{S})\setminus F(S,\widetilde{S})|};~\tau=+\infty\Big)
&\leq \big(nC_1\big)^{\lfloor\theta\log n\rfloor}\big(\frac{1}{n(1-\delta_2)}\big)^{\lfloor\theta \log n\rfloor-1} \\
&=n\big(\frac{C_1}{1-\delta_2}\big)^{\lfloor\theta\log n\rfloor}<1 \notag
\end{align}
for sufficiently large $n$, since $\big(\frac{C_1}{1-\delta_2}\big)^\theta<\frac{1}{e}$ as we have introduced. By Equations \eqref{equ 4.18} and \eqref{equ 4.19},
\begin{equation}\label{equ 4.20}
\widetilde{\rm E}\Big(\big(nC_1\big)^{|F(S,\widetilde{S})|}~C_2^{|D(S,\widetilde{S})\setminus F(S,\widetilde{S})|}\Big)<1+\frac{C_2}{1-\gamma}\big(\frac{1}{1-\frac{C_1}{1-\delta_2}}\big)^2<+\infty
\end{equation}
for sufficiently large $n$.  Equation \eqref{equ 4.8} follows from Equation \eqref{equ 4.20} directly and the proof is complete.

\qed

\quad

\textbf{Acknowledgments.} The author is grateful to the financial
support from the National Natural Science Foundation of China with
grant number 11501542.

{}
\end{document}